\documentclass[12pt, reqno]{amsart}
\usepackage{mathrsfs}
\renewcommand{\mathcal}{\mathscr}
\usepackage{amsmath,amsthm,amsopn,amstext,amscd,amsfonts,amssymb}
\usepackage{graphicx, color, enumerate}
\usepackage[latin1]{inputenc}
\usepackage[active]{srcltx}
\usepackage{tikz}
\usepackage{pgf}
\usepackage{etex}
\usepackage{verbatim}
\usepackage{tikz-3dplot}
\usepackage{pgfkeys}
\usepackage{fp}
\usepackage[shortlabels]{enumitem}
\usepackage[colorlinks=true,linkcolor=blue]{hyperref}
\numberwithin{equation}{section}

\numberwithin{bbb}{section}

\newtheorem{theorem}{Theorem}
\newtheorem{definition}{Definition}
\newtheorem{lemma}{Lemma}

\newtheorem{remark}{Remark}

\title[Global existence band Asymptotic behavior]{Global existence and asymptotic behavior for a reaction-diffusion system with unbounded coefficients}
\author[M. Majdoub, N. Tatar]{Mohamed Majdoub and Nasser-eddine Tatar}
\address[M. Majdoub]{Department of Mathematics, College of Science, Imam Abdulrahman Bin Faisal University, P. O. Box 1982, Dammam, Saudi Arabia. \newline Basic and Applied Scientific Research Center, Imam Abdulrahman Bin Faisal University, P.O. Box 1982, 31441, Dammam, Saudi Arabia.}
\email{\sl mmajdoub@iau.edu.sa}
\address[N. Tatar]{King Fahd University of Petroleum and Minerals, Department of Mathematics, IRC for Intelligent Manufacturing and Robotics, Dhahran, 31261, Saudi Arabia.}
\email{\sl tatarn@kfupm.edu.sa}
\begin{document}

\begin{abstract}
We consider a reaction-diffusion system which may serve as a model for a ferment catalytic
reaction in chemistry. The model consists of a system of reaction diffusion
equations with unbounded time dependent coefficients and different
polynomial reaction terms. An exponential decay of the globally bounded
solutions is proved. The key tool of the proofs are properties of analytic
semigroups and some inequalities.
\end{abstract}

\subjclass[2010]{35K57, 35A01, 35K57,35B40, 35Q92}

\keywords{Analytic semi-group, exponential decay, fractional
operator, reaction-diffusion system, sectorial operator.}

\date{\today}

\maketitle

\section{Introduction}
In the past three decades a strong interest in the questions of global existence and asymptotic behavior of solutions to various classes of reaction-diffusion systems has been witnessed. A particular attention was paid to autonomous systems and different sufficient conditions were assumed to ensure global existence of classical solutions. A main difficulty in proving global existence is the lack of maximum principle estimates or invariant regions. See, among many, \cite{Am1985, 2, 1, LP2017,
8, 10, PS1997, 12}. Recently, the global existence of classical solutions was shown in \cite{FMT2020} for quasi-positive nonlinearities having a (slightly super-) quadratic growth and obeying a mass control assumption. This, in particular, includes the cases of mass conservation and mass dissipation. See also \cite{LP2020} for initial data of low regularity.

Our goal in this paper is to study the following non-autonomous system
\begin{align} \label{1}
	\left\{
	\begin{array}{ccl}
u_{t}-(d_{1}\Delta -b_{1})u&=&a_{1}(t)w^{m}-a_{2}(t)u^{n}v^{k}, \\
v_{t}-d_{2}\Delta v&=&\left( a_{1}(t)+a_{3}(t)\right) w^{m}-a_{2}(t)u^{n}v^{k},
\\
w_{t}-(d_{3}\Delta -b_{2})w&=&-\left( a_{1}(t)+a_{3}(t)\right)
w^{m}-a_{4}(t)w+a_{2}(t)u^{n}v^{k},  \\
\frac{\partial u}{\partial \nu }=\frac{\partial v}{\partial \nu }=\frac{%
\partial w}{\partial \nu }&=&0, \\
(u,v,w)(x,0)&=&(u_{0}(x),v_{0}(x),w_{0}(x)),
\end{array}
	\right.
	\end{align}
where $t>0$, $x\in\Omega \subset \mathbf{R}^{N},$ $N\geq 1$, $\Omega$ is a bounded domain with
smooth boundary $\partial \Omega $ and $\partial /\partial \nu $ denotes the
outward normal derivative on $\partial \Omega .$ The diffusion coefficients $%
d_{i},$ $i=1,2,3$ are assumed positive constants and $u_{0},$ $v_{0}$ and $%
w_{0}$ are given nonnegative bounded initial functions. The coefficients $%
a_{i}(t)$ are of the form $t^{\sigma _{i}}h_{i}(t),$ where $\sigma _{i}\geq
0,$ $i=1,2,3$ and $\sigma _{4}=0.$ The functions $h_{i}(t)$ are nonnegative
continuous functions. The powers $m,$ $n$ and $k$ are assumed greater than
one.

Systems with time dependent nonlinearities were considered in \cite{5} by Kahane.
Specifically, a system of the form
\begin{equation*}
\left\{
\begin{array}{l}
-u_{t}+Lu=f(x,t,u,v)\;in\;\Omega \times (0,\infty ) \\
-v_{t}+Mv=g(x,t,u,v)\;in\;\Omega \times (0,\infty )%
\end{array}
\right.
\end{equation*}
where $L$ and $M$ are uniformly elliptic operators, with boundary conditions
of Robin type is studied. He proved that the solution converges to the
stationary state, that is the solution of the limiting elliptic problem. To
this end he assumed that
\begin{equation*}
f(x,t,u,v)\rightarrow \bar{f}(x,u,v)
\end{equation*}
and
\begin{equation*}
g(x,t,u,v)\rightarrow \bar{g}(x,u,v)
\end{equation*}
uniformly in $\Omega $ and $(u,v)$ in any bounded subset of the first
quadrant in $\mathbf{R}^{2}.$ The matrix formed by the partial derivatives $%
\bar{f}_{u},\bar{f}_{v},\bar{g}_{u}$ and $\bar{g}_{v}$ satisfies a column
diagonal dominance type condition. In the present work we do not make such
restrictions.

Let us return now to our system \eqref{1}. In the case $m=n=k=1,$ $b_{1}=b_{2}=0$ and $a_{i}(t)\equiv a_{i}$ are constants, this problem has been studied by Wang \cite{14}. If moreover
$a_{3}=a_{4}=0,$ then this problem may also serve as a model for sugar transporting into red blood cells (see Ruan \cite{13}, Rothe \cite{12}, Feng \cite{1}, Morgan \cite{10} and references therein). In particular, Wang proved a convergence result and an exponential decay result in $C^{\mu }(\bar{\Omega})$, $\mu \in [0,2).$ Unfortunately, his method seems to be not valid for
the present problem \eqref{1} because of the unboundedness of the coefficients.

To our knowledge, a system with time-dependent coefficients like \eqref{1} is still not well studied. Nevertheless,
using similar ideas as in \cite{FMT2020}, we are able to prove global existence of classical solutions to \eqref{1} provided that $0\leq m\leq 2+\varepsilon$, $n, k \geq 0$ and $n+k\leq 2+\varepsilon$ for sufficiently small $\varepsilon>0$. More precisely, we assume that
\begin{equation}
\label{GE1}
0\leq m\leq 2+\varepsilon\quad\mbox{and}\quad n+k\leq 2+\varepsilon,
\end{equation}
\begin{equation}
\label{GE2}
u_0, v_0, w_0 \in L^1(\Omega)\cap L^\infty(\Omega)\quad\mbox{and}\quad u_0, v_0, w_0\geq 0.
\end{equation}
\begin{theorem}
\label{GE}
Suppose that assumption \eqref{GE2} is fulfilled. Then, there exists $\varepsilon>0$ such that if \eqref{GE1} is satisfied, \eqref{1} possesses a unique classical global solution.
\end{theorem}

Our main task here is to study the asymptotic behavior and the decay rate of global classical solutions obtained in Theorem \ref{GE}. To this end we shall adapt the methods used in Hoshino \cite{3} for the former task and the methods used in Kirane and Tatar \cite{7} for the
latter goal. To state our results in a clear way, we shall fix some notations. For $\alpha\in (0,1)$, we denote by
\begin{equation*}
q^{*}=q^{*}(\alpha)=\left\{
\begin{array}{l}
\frac{2}{\alpha}-1,\;if\;\frac{1}{2}\leq \alpha<1 \\
2,\;if\;0<\alpha <\frac{1}{2}%
\end{array}
\right.
\end{equation*}
and $q$ the Lebesgue conjugate exponent of $q^{*}$. Let $r=\frac{1-\alpha}{\alpha}$.
First, we state a convergence result of the solutions in the space $C^{\mu }(\Omega ),$ $\mu \in [0,2).$ This result will be needed in the sequel. Without loss of generality, we shall assume that $b_{1}=b_{2}=b.$
\begin{theorem}
\label{T1}
Suppose that $1-q\alpha >0$, $1+q(\sigma _{i}-\alpha l)>0,$ $i=1,2,3$ for
some $l$ such that $1<l<\min \{m,n,k\},$ $h_{i}\in L^{q^{*}}(0,\infty ),$ $%
i=1,2,3,4$ and $a_{1}(t)\leq Ca_{3}(t)$ for some positive constant $C$, $%
\forall t>0.$ Then, for every $\mu \in [0,2)$%
\begin{equation*}
u(t)\rightarrow 0,v(t)\rightarrow v_{\infty }\;and\;w(t)\rightarrow
0\;in\;C^{\mu }(\Omega )\;as\;t\rightarrow \infty
\end{equation*}
where $v_{\infty }=\left| \Omega \right| ^{-1}\left\{ \int\limits_{\Omega
}w_{0}dx-\int\limits_{0}^{\infty }\int\limits_{\Omega
}(h_{4}(s)+b)w(s)dxds\right\} .$
\end{theorem}
\begin{remark}
{\rm The proof of Theorem \ref{T1} is similar to the proof of Theorem 2.2 in \cite{3} with minor
modifications that will be clear from the proof of our next and main result.
It is therefore omitted.}
\end{remark}

To state our next result, we suppose that $N$ and $p$ satisfy
\begin{equation}
\label{H}
2\alpha >\frac{N}{p}\quad\mbox{and}\quad \max \left\{ 1,\frac{N(m-1)}{2p\alpha },\frac{N(n-1)}{2p\alpha }\right\} <\min \{m,n\}.
\end{equation}

\begin{theorem}
\label{T2}
Assume that the hypotheses of Theorem \ref{T1} and \eqref{H} are fulfilled. If
\begin{equation*}
(l-1)\int\limits_{0}^{\infty }h(s)ds<\log \left( \frac{1+C_{0}^{l-1}}{C_{0}}%
\right)
\end{equation*}
for some constant $l$ such that $\max \left\{ 1,\frac{N(m-1)}{2p\alpha },
\frac{N(n-1)}{2p\alpha }\right\} <l<\min \{m,n\}$ and a constant $%
C_{0}=C_{0}(\left\| u_{0}\right\| _{p},\left\| w_{0}\right\| _{p})$ and a function $h$ to be
determined in the proof, then for $\mu \in [0,2)$ and $N,$ $p$ such that $
0\leq \mu <2\alpha -\frac{N}{p}$ we have
\begin{enumerate}
\item[(a)] $\left\| u\right\| _{C^{\mu }(\bar{\Omega})},$ $\left\| w\right\|
_{C^{\mu }(\bar{\Omega})}\leq Ce^{-(b-\varepsilon )t}\left( \left\|
u_{0}\right\| _{p}+\left\| w_{0}\right\| _{p}\right) $, where $0<\varepsilon
<b$, as $t\rightarrow \infty .$

\item[(b)]
\begin{itemize}
\item[(i)] If $d_{2}\lambda <lb,$%
\begin{equation*}
\left\| v-v_{\infty }\right\| _{C^{\mu }(\bar{\Omega})}\leq Ce^{-\min
\left\{ (b-\varepsilon ),d_{2}\lambda \right\} t}\;as\;t\rightarrow \infty .
\end{equation*}
\item[(ii)] If $d_{2}\lambda \geq lb$ and for $i=1,2,3$
\begin{equation*}
\int\limits_{0}^{t}e^{q^{*}\rho s}h_{i}^{q^{*}}(s)ds=O(e^{q^{*}\tilde{\rho}%
t})\;as\;t\rightarrow \infty
\end{equation*}
for some $\rho >d_{2}\lambda -l(b-\varepsilon )$ and $\tilde{\rho}%
<d_{2}\lambda ,$ then
\begin{equation*}
\left\| v-v_{\infty }\right\| _{C^{\mu }(\bar{\Omega})}\leq Ce^{-\min
\left\{ (b-\varepsilon ),d_{2}\lambda -\tilde{\rho}\right\}
t}\;as\;t\rightarrow \infty .
\end{equation*}
\end{itemize}
\end{enumerate}
\end{theorem}

\section{Preliminaries}

In this section we present the notation that will be used in this paper and
prepare some material which will be useful in our proofs.

By $W^{l,p}(\Omega )$ we denote the usual Sobolev space of order $l\geq 0$
for $1\leq p\leq \infty $. The space $C^{\mu }(\bar{\Omega}),\mu \geq 0$ is
the Banach space of $[\mu ]-$times continuously differentiable functions in $%
\bar{\Omega}$ whose $[\mu ]-$th order derivatives are H\"{o}lder continuous
with exponent $\mu -[\mu ].$

\begin{definition}
For $p\in (1,\infty )$, we define

$D(A_{p})=D(B_{p})=D(G_{p})=\left\{ y\in W^{2,p}(\Omega ):\partial
y/\partial \nu \left| _{\partial \Omega }\right. =0\right\} ,$

$A_{p}y=-(d_{1}\Delta -b_{1})y,$ $B_{p}y=-d_{2}\Delta y,$ $%
G_{p}y=-(d_{3}\Delta -b_{2})y.$
\end{definition}

The operators $-A_{p},$ $-B_{p}$ and $-G_{p}$ defined in this way are
sectorial operators (see Henry \cite{2}) and generate analytic semigroups
$\left\{ e^{-tA_{p}}\right\} _{t\geq 0},$ $\left\{ e^{-tB_{p}}\right\}
_{t\geq 0}$ and $\left\{ e^{-tG_{p}}\right\} _{t\geq 0},$ respectively.

\begin{definition}
By $Q_{0}$ and $Q_{+}$ we designate the following projection operators, for $%
y\in L^{p}(\Omega )$, $p\in (1,\infty )$%
\begin{equation*}
Q_{0}y=\left| \Omega \right| ^{-1}\int\limits_{\Omega
}y(x)dx\;and\;Q_{+}y=y-Q_{0}y
\end{equation*}
where $\left| \Omega \right| $ is the volume of $\Omega .$
\end{definition}

\begin{definition}
The restriction of $B_{p}$ onto $Q_{+}L^{p}(\Omega )$ will be denoted by $%
B_{p+}$ i.e. $B_{p+}=B_{p}\left| _{Q_{+}L^{p}(\Omega )}\right. .$
\end{definition}

For $p\in (1,\infty ),$ the operator $B_{p+}$ generates an analytic
semigroup denoted by $\left\{ e^{-tB_{p+}}\right\} _{t\geq 0}$ in $%
Q_{+}L^{p}(\Omega ).$

We define the fractional powers of the above operators in the usual way (see
Henry \cite{2}).

\begin{lemma}
\label{L1}
Let $A$ be a sectorial operator in $X=L^{p}(\Omega ),$ $1\leq p<\infty $
with $D(A)=X^{1}\subset W^{m,p}(\Omega )$ for some $m\geq 1.$ Then, for $%
0\leq \alpha \leq 1$, we have $X^{\alpha }\subset C^{\mu }$ when $0\leq \mu
<m\alpha -\frac{N}{p}.$
\end{lemma}

\begin{lemma}
\label{L2}
Let $\lambda $ denote the least positive eigenvalue of the Laplacian with
homogeneous Neumann boundary condition. Let $p\in (1,\infty )$. For every $%
\alpha \in [0,1)$ there exist positive constants $C_{i},$ $i=1,2,3$ such
that for $t>0$ and $y\in L^{p}(\Omega )$

$\left\| A_{p}^{\alpha }e^{-tA_{p}}y\right\| _{p}\leq C_{1}t^{-\alpha
}e^{-b_{1}t}\left\| y\right\| _{p},$

$\left\| B_{p+}^{\alpha }e^{-tB_{p+}}Q_{+}y\right\| _{p}\leq C_{2}t^{-\alpha
}e^{-d_{2}\lambda t}\left\| Q_{+}y\right\| _{p},$

$\left\| G_{p}^{\alpha }e^{-tG_{p}}y\right\| _{p}\leq C_{3}t^{-\alpha
}e^{-b_{2}t}\left\| y\right\| _{p}.$
\end{lemma}

\begin{lemma}
\label{L3}
Let $p\in (1,\infty ),$ $l\geq 1$ and $\alpha \in [0,1).$ Then, there exists
a positive constant $C$ such that
\begin{equation*}
\left\| y\right\| _{pl}\leq C\left\| A_{p}^{\alpha }y\right\| _{p}^{\theta
}.\left\| y\right\| _{p}^{1-\theta }
\end{equation*}
where $\theta $ satisfies
\begin{equation*}
\frac{N(l-1)}{2pl\alpha }<\theta <1.
\end{equation*}
\end{lemma}

See Henry \cite{2} for the proofs of Lemmas \ref{L1}-\ref{L2}-\ref{L3}.

\begin{lemma}
\label{L4}
Let $\alpha \in [0,1)$ and $\beta \in \mathbf{R}.$ There exists a positive
constant $C=C(\alpha ,\beta )$ such that
\begin{equation*}
\int\limits_{0}^{t}s^{-\alpha }e^{\beta s}ds\leq \left\{
\begin{array}{l}
Ce^{\beta t}\;if\;\beta >0 \\
C(t+1)\;if\;\beta =0 \\
C\;if\;\beta <0.%
\end{array}
\right.
\end{equation*}
\end{lemma}

See Hoshino and Yamada \cite{4}, for instance, for the proof of Lemma \ref{L4}.

\begin{lemma}
\label{L5}
Let $\mu ,$ $\nu ,$ $\tau $ and $z>0$, then
\begin{equation*}
z^{1-\nu }\int\limits_{0}^{z}(z-\xi )^{\nu -1}\xi ^{\mu -1}e^{-\tau \xi
}d\xi \leq C\tau ^{-\mu }
\end{equation*}
where $C$ is a constant independent of $z.$
\end{lemma}

See Michalski \cite{9} or Kirane and Tatar \cite{6} for the proof of Lemma \ref{L5}.\\

Let $y : [a,b]\to [0,\infty)$  be a continuous function satisfying
\begin{equation}
\label{GI1}
y(t)\leq M+\int\limits_{a}^{t}\,\lambda(s)\,g(y(s))ds,\;a\leq t\leq b
\end{equation}
where $M>0$, $\lambda : [a,b]\to [0,\infty)$ is continuous and $g : [0,\infty)\to [0,\infty)$ is continuous nondecreasing with
$$
g(y)>0,\quad \forall\;\;\;y>0.
$$
We define
\begin{equation}
\label{GM}
G_{M}(y)=\int\limits_{M}^{y}\frac{ds}{g(s)},\;y>0.
\end{equation}

The following Bihari-type inequality is classical and can be found in many references. See for instance \cite{11}.
\begin{lemma}
\label{L6}
Under the above assumptions on $y, \lambda,  g$, we have
\begin{equation}
\label{Bihari}
y(t)\leq G_M^{-1}\bigg(\int\limits_a^t\,\lambda(s)\,ds\bigg),
\end{equation}
for $a\leq t\leq T\leq b$ with $T>a$ satisfies
\begin{equation}
\label{T}
\int\limits_a^T\,\lambda(s)\,ds<\int\limits_M^\infty\,\frac{ds}{g(s)}\,ds.
\end{equation}
\end{lemma}
\begin{remark}
{\rm The classical Gronwall's inequality corresponds to $g(y)=y$.}
\end{remark}

\section{Proof of Theorem \ref{GE}}
Define $\mathbf{z}=(u,v,w)$, we rewrite \eqref{1} as
\begin{align} \label{z-y}
	\left\{
	\begin{array}{ccl}
	u_{t}-d_{1}\Delta\,u&=&f_1(t,\mathbf{z}), \\
v_{t}-d_{2}\Delta v&=&f_2(t,\mathbf{z}),\\
w_{t}-d_{3}\Delta\,w&=&f_3(t,\mathbf{z}), \\
y_t-\Delta y&=&f_4(t,\mathbf{z}),\;\; y(0)=0,
	\end{array}
	\right.
	\end{align}
where $f_4:= -f_1-f_2-f_3$ and the last equation in $y$ is added to fulfill technical requirements as it will be clear below. From \eqref{z-y} we deduce that
\begin{equation}
\label{uvwy}
 \left(u+v+w+y\right)_t-\Delta\left(d_1 u+d_2 v+d_3 w+y\right)=0.
\end{equation}
It follows that
\begin{equation}
 \label{uvwy1}
 u(x,t)+v(x,t)+w(x,t)+y(x,t)-\Delta \mathbf{m}(x,t)=g(x),
\end{equation}
where
\begin{equation}
\label{m}
\mathbf{m}(x,t)=\int_0^t\,\left[d_1 u(x,s)+d_2 v(x,s)+d_3 w(x,s)+y(x,s)\right]\,ds,
\end{equation}
and
\begin{equation}
\label{g}
g(x)=u_0(x)+v_0(x)+w_0(x).
\end{equation}
Therefore
\begin{equation}
\label{Est-Linfty}
\|\mathbf{z}\|_{L^\infty_{x,t}}\leq \|\Delta\mathbf{m}\|_{L^\infty_{x,t}}+\|g\|_{L^\infty_{x}}\quad\text{and}\quad \|\mathbf{m}\|_{L^\infty_{x,t}}\leq \mathcal{C}_T,
\end{equation}
where $\mathcal{C}_T$ depends continuously in $T>0$.
Using \eqref{uvwy1}, one can write
\begin{equation}
\label{Eq-m}
B(x,t)\partial_t\mathbf{m}-\Delta\mathbf{m}=g(x),
\end{equation}
where
$$
0<\underline{B}\leq B(x,t)=\frac{u(x,t)+v(x,t)+w(x,t)+y(x,t)}{d_1 u(x,t)+d_2 v(x,t)+d_3 w(x,t)+y(x,t)}\leq \overline{B}<\infty.
$$
Lemma 3.5 in \cite{FMT2020} implies
\begin{equation}
\label{Holder-m}
|\mathbf{m}(x,t)-\mathbf{m}(x',t)|\leq \mathcal{C}_T |x-x'|^{\delta},\;\;\; \forall\; 0<t<T,
\end{equation}
for some $\delta\in (0,1)$. Equation \eqref{uvwy1} also implies
\begin{equation}
\label{Eq-mm}
\partial_t\mathbf{m}-\Delta\mathbf{m}=(d_1-1) u+(d_2-1) v+(d_3-1) w+g(x).
\end{equation}
By using the H\"older continuity of $\mathbf{m}$ and \cite[Lemma 1.1, p. 283]{FMT2020}, we infer
\begin{equation}
\label{Est-mm}
|\nabla\mathbf{m}|\leq \mathcal{C}_T\,|\mathbf{z}|^{\frac{1-\delta}{2-\delta}}.
\end{equation}
Applying \cite[Lemma 1.1, p. 283]{FMT2020} to each equation of $q\in \{u,v,w,y\}$ leads to
\begin{equation}
\label{Est-q}
|\nabla\,q|\leq \mathcal{C}_T\,\left(1+|\mathbf{z}|\right)^{\frac{3}{2}}.
\end{equation}
Using \cite[Lemma 3.9, p. 295]{FMT2020} together with \eqref{Est-mm} and \eqref{Est-q} gives
\begin{eqnarray*}
|\Delta\mathbf{m}|&\leq& \mathcal{C}_T\, |\nabla\mathbf{m}|^{1/2}\,\left(|\nabla\,u|+ |\nabla\,v|+ |\nabla\,w|+\, |\nabla\, g|\right)^{1/2}\\
&\leq&\mathcal{C}_T\,|\mathbf{z}|^{\frac{1-\delta}{2(2-\delta)}}\,\left(1+|\mathbf{z}|^{3/2}\right)^{1/2}\\
&\leq&\mathcal{C}_T\,\left(1+|\mathbf{z}|^{\frac{1-\delta}{2(2-\delta)}+\frac{3}{4}}\right).
\end{eqnarray*}
From the above estimate and \eqref{Est-Linfty}, we obtain that
$$
\|\mathbf{z}\|_{L^\infty_{x,t}}\leq \mathcal{C}_T\,\left(1+\|\mathbf{z}\|_{L^\infty_{x,t}}^{\frac{1-\delta}{2(2-\delta)}+\frac{3}{4}}\right).
$$
Since $\frac{1-\delta}{2(2-\delta)}+\frac{3}{4}<1$ then
$\|\mathbf{z}\|_{L^\infty_{x,t}}\leq \mathcal{C}_T$. Consequently $T_{max}=\infty$.
\section{Proof of Theorem \ref{T2}}

\begin{proof}

\underline{A. The decay rate of $\left\| u\right\| _{C^{\mu }(\bar{\Omega})}$
and $\left\| w\right\| _{C^{\mu }(\bar{\Omega})}$}:

Clearly we have the integral equations associated with the first and third
equation of \eqref{1}
\begin{equation}
u(t)=e^{-tA_{p}}u(0)+\int\limits_{0}^{t}e^{-(t-s)A_{p}}\left\{
a_{1}(s)w^{m}-a_{2}(s)u^{n}v^{k}\right\} ds  \label{3}
\end{equation}
\begin{equation}
\label{4}
\begin{array}{c}
w(t)=e^{-tG_{p}}w(0) \\
+\int\limits_{0}^{t}e^{-(t-s)G_{p}}\left\{
-(a_{1}(s)+a_{3}(s))w^{m}-a_{4}(s)w+a_{2}(s)u^{n}v^{k}\right\} ds.%
\end{array}%
\end{equation}
Applying $A_{p}^{\alpha }$ and $G_{p}^{\alpha }$, $0<\alpha <1$ to both
sides of \eqref{3} and \eqref{4} respectively, we infer from Lemma \ref{L5} that
\begin{equation*}
\begin{array}{c}
\left\| A_{p}^{\alpha }u\right\| _{p}\leq C_{1}t^{-\alpha }e^{-bt}\left\|
u_{0}\right\| _{p}+C_{1}\int\limits_{0}^{t}(t-s)^{-\alpha
}e^{-b(t-s)}s^{\sigma _{1}}h_{1}(s)\left\| w^{m}\right\| _{p}ds \\
+C_{1}\int\limits_{0}^{t}(t-s)^{-\alpha }e^{-b(t-s)}s^{\sigma
_{2}}h_{2}(s)\left\| u^{n}v^{k}\right\| _{p}ds%
\end{array}
\end{equation*}
and
\begin{equation*}
\begin{array}{ll}
\left\| G_{p}^{\alpha }w\right\| _{p}\leq & C_{3}t^{-\alpha }e^{-bt}\left\|
w_{0}\right\| _{p} \\
& +C_{3}\int\limits_{0}^{t}(t-s)^{-\alpha }e^{-b(t-s)}(s^{\sigma
_{1}}h_{1}(s)+s^{\sigma _{3}}h_{3}(s))\left\| w^{m}\right\| _{p}ds \\
& +C_{3}\int\limits_{0}^{t}(t-s)^{-\alpha }e^{-b(t-s)}h_{4}(s)\left\|
w\right\| _{p}ds \\
& +C_{3}\int\limits_{0}^{t}(t-s)^{-\alpha }e^{-b(t-s)}a_{2}(s)\left\|
u^{n}v^{k}\right\| _{p}ds.%
\end{array}
\end{equation*}
As $0\leq v\leq M,$ for some $M>0,$ we have for all $t\geq \delta >0$%
\begin{equation}
\label{5}
\begin{array}{c}
e^{bt}\left\| A_{p}^{\alpha }u\right\| _{p}\leq C_{1}\delta ^{-\alpha
}\left\| u_{0}\right\| _{p}+C_{1}\int\limits_{0}^{t}(t-s)^{-\alpha
}e^{bs}s^{\sigma _{1}}h_{1}(s)\left\| w^{m}\right\| _{p}ds \\
+C_{1}M^{k}\int\limits_{0}^{t}(t-s)^{-\alpha }e^{bs}s^{\sigma
_{2}}h_{2}(s)\left\| u^{n}\right\| _{p}ds,%
\end{array}%
\end{equation}
and
\begin{equation}
\label{6}
\begin{array}{c}
e^{bt}\left\| G_{p}^{\alpha }w\right\| _{p}\leq C_{3}\delta ^{-\alpha
}\left\| w_{0}\right\| _{p}+C_{3}\int\limits_{0}^{t}(t-s)^{-\alpha
}e^{bs}h_{4}(s)\left\| w\right\| _{p}ds \\
+C_{3}\int\limits_{0}^{t}(t-s)^{-\alpha }e^{bs}(s^{\sigma
_{1}}h_{1}+s^{\sigma _{3}}h_{3})(s)\left\| w^{m}\right\| _{p}ds \\
+C_{3}M^{k}\int\limits_{0}^{t}(t-s)^{-\alpha }e^{bs}s^{\sigma
_{2}}h_{2}(s)\left\| u^{n}\right\| _{p}ds.%
\end{array}%
\end{equation}
Lemma \ref{L6} and the uniform boundedness of $w$ allow us to write
\begin{equation}
\left\| w^{m}\right\| _{p}=\left\| w\right\| _{mp}^{m}\leq C\left\|
G_{p}^{\alpha }w\right\| _{p}^{m\theta }.\left\| w\right\| _{p}^{m(1-\theta
)}\leq C_{4}\left\| G_{p}^{\alpha }w\right\| _{p}^{m\theta }  \label{7}
\end{equation}
with $\frac{N(m-1)}{2pm\alpha }<\theta <1.$

In the same way we have
\begin{equation}
\left\| u^{n}\right\| _{p}\leq C_{5}\left\| A_{p}^{\alpha }u\right\|
_{p}^{n\theta }\;with\;\frac{N(n-1)}{2pn\alpha }<\theta <1.  \label{8}
\end{equation}
Let us choose $l$ such that
\begin{equation*}
\max \left\{ 1,\frac{N(m-1)}{2p\alpha },\frac{N(n-1)}{2p\alpha }\right\}
<l<\min \{m,n\}
\end{equation*}
and $\theta =\frac{l}{m}$ in \eqref{7} and $\theta =\frac{l}{n}$ in \eqref{8}. Then
\begin{equation}
\left\| w^{m}\right\| _{p}\leq C_{4}\left\| G_{p}^{\alpha }w\right\|
_{p}^{l}\;and\;\left\| u^{n}\right\| _{p}\leq C_{5}\left\| A_{p}^{\alpha
}u\right\| _{p}^{l}.  \label{9}
\end{equation}
Taking \eqref{9} into account in \eqref{5} and \eqref{6}, we obtain
\begin{equation}
\begin{array}{c}
e^{bt}\left\| A_{p}^{\alpha }u\right\| _{p}\leq C_{1}\delta ^{-\alpha
}\left\| u_{0}\right\| _{p} \\
+C_{6}\int\limits_{0}^{t}(t-s)^{-\alpha }e^{b(1-l)s}s^{\sigma
_{1}}h_{1}(s)\left( e^{bs}\left\| G_{p}^{\alpha }w\right\| _{p}\right) ^{l}ds
\\
+C_{7}\int\limits_{0}^{t}(t-s)^{-\alpha }e^{b(1-l)s}s^{\sigma
_{2}}h_{2}(s)\left( e^{bs}\left\| A_{p}^{\alpha }u\right\| _{p}\right)
^{l}ds.%
\end{array}
\label{10}
\end{equation}
and
\begin{equation}
\begin{array}{c}
e^{bt}\left\| G_{p}^{\alpha }w\right\| _{p}\leq C_{3}\delta ^{-\alpha
}\left\| w_{0}\right\| _{p}+C_{3}\int\limits_{0}^{t}(t-s)^{-\alpha
}h_{4}(s)e^{bs}\left\| w\right\| _{p}ds \\
+C_{8}\int\limits_{0}^{t}(t-s)^{-\alpha }e^{b(1-l)s}(s^{\sigma
_{1}}h_{1}+s^{\sigma _{3}}h_{3})(s)\left( e^{bs}\left\| G_{p}^{\alpha
}w\right\| _{p}\right) ^{l}ds \\
+C_{9}\int\limits_{0}^{t}(t-s)^{-\alpha }e^{b(1-l)s}s^{\sigma
_{2}}h_{2}(s)\left( e^{bs}\left\| A_{p}^{\alpha }u\right\| _{p}\right)
^{l}ds.%
\end{array}
\label{11}
\end{equation}
The second term in the right hand side of \eqref{11} may be estimated in the
following manner, for $0<\varepsilon <b$
\begin{equation}
\begin{array}{c}
\int\limits_{0}^{t}(t-s)^{-\alpha }h_{4}(s)e^{bs}\left\| w\right\|
_{p}ds=\int\limits_{0}^{t}(t-s)^{-\alpha }e^{\varepsilon s}e^{(b-\varepsilon
)s}h_{4}(s)\left\| w\right\| _{p}ds \\
\leq \left( \int\limits_{0}^{t}(t-s)^{-\alpha q}e^{\varepsilon qs}ds\right)
^{\frac{1}{q}}\left( \int\limits_{0}^{t}h_{4}^{q^{*}}(s)\left(
e^{(b-\varepsilon )s}\left\| w\right\| _{p}\right) ^{q^{*}}ds\right) ^{\frac{%
1}{q^{*}}} \\
\leq C_{10}e^{\varepsilon t}\left( \int\limits_{0}^{t}h_{4}^{q^{*}}(s)\left(
e^{(b-\varepsilon )s}\left\| w\right\| _{p}\right) ^{q^{*}}ds\right) ^{\frac{%
1}{q^{*}}}.%
\end{array}
\label{12}
\end{equation}
We have used the H\"{o}lder inequality, Lemma 7 and the embedding $%
D(G_{p}^{\alpha })\subset L^{p}$ to derive the last inequalities in \eqref{12}.
Multiplying \eqref{10} and \eqref{11} by $e^{-\varepsilon t},$ setting
\begin{equation*}
U(t)=e^{(b-\varepsilon )t}\left\| A_{p}^{\alpha }u\right\| _{p},\text{ }%
W(t)=e^{(b-\varepsilon )t}\left\| G_{p}^{\alpha }w\right\| _{p}
\end{equation*}
and taking into account \eqref{12} we find for $t\geq \delta >0$%
\begin{equation*}
\begin{array}{c}
U(t)\leq C_{1}\delta ^{-\alpha }\left\| u_{0}\right\|
_{p}+C_{6}\int\limits_{0}^{t}(t-s)^{-\alpha }e^{-(b-\varepsilon
)(l-1)s}s^{\sigma _{1}}h_{1}(s)W(s)^{l}ds \\
+C_{7}\int\limits_{0}^{t}(t-s)^{-\alpha }e^{-(b-\varepsilon
)(l-1)s}s^{\sigma _{2}}h_{2}(s)U(s)^{l}ds,%
\end{array}
\end{equation*}
and
\begin{equation*}
\begin{array}{c}
W(t)\leq C_{3}\delta ^{-\alpha }\left\| w_{0}\right\| _{p}+C_{10}\left(
\int\limits_{0}^{t}h_{4}^{q^{*}}(s)W(s)^{q^{*}}ds\right) ^{\frac{1}{q^{*}}}
\\
+C_{8}\int\limits_{0}^{t}(t-s)^{-\alpha }e^{-(b-\varepsilon
)(l-1)s}(s^{\sigma _{1}}h_{1}+s^{\sigma _{3}}h_{3})(s)W(s)^{l}ds \\
+C_{9}\int\limits_{0}^{t}(t-s)^{-\alpha }e^{-(b-\varepsilon
)(l-1)s}s^{\sigma _{2}}h_{2}(s)U(s)^{l}ds.%
\end{array}
\end{equation*}

Next, as $1-q\alpha >0,$ then using the H\"{o}lder inequality and Lemma 8,
we obtain for $t\geq \delta >0$%
\begin{equation*}
\begin{array}{c}
U(t)\leq C_{1}\delta ^{-\alpha }\left\| u_{0}\right\| _{p}+C_{11}\delta
^{-\alpha }\left( \int\limits_{0}^{t}h_{1}^{q^{*}}(s)W(s)^{lq^{*}}ds\right)
^{\frac{1}{q^{*}}} \\
+C_{12}\delta ^{-\alpha }\left(
\int\limits_{0}^{t}h_{2}^{q^{*}}(s)U(s)^{lq^{*}}ds\right) ^{\frac{1}{q^{*}}}%
\end{array}
\end{equation*}
and
\begin{equation*}
\begin{array}{c}
W(t)\leq C_{3}\delta ^{-\alpha }\left\| w_{0}\right\| _{p}+C_{10}\left(
\int\limits_{0}^{t}h_{4}^{q^{*}}(s)W(s)^{q^{*}}ds\right) ^{\frac{1}{q^{*}}}
\\
+C_{13}\delta ^{-\alpha }\left(
\int\limits_{0}^{t}h_{1}^{q^{*}}(s)W(s)^{lq^{*}}ds\right) ^{\frac{1}{q^{*}}%
}+C_{14}\delta ^{-\alpha }\left(
\int\limits_{0}^{t}h_{3}^{q^{*}}(s)W(s)^{lq^{*}}ds\right) ^{\frac{1}{q^{*}}}
\\
+C_{15}\delta ^{-\alpha }\left(
\int\limits_{0}^{t}h_{2}^{q^{*}}(s)U(s)^{lq^{*}}ds\right) ^{\frac{1}{q^{*}}}.%
\end{array}
\end{equation*}
Now we use the inequality
\begin{equation*}
(x_{1}+x_{2}+...+x_{n})^{r}\leq n^{r-1}(x_{1}^{r}+x_{2}^{r}+...+x_{n}^{r})
\end{equation*}
to get
\begin{equation}
\begin{array}{c}
U(t)^{q^{*}}\leq 3^{q^{*}-1}\left( C_{1}\delta ^{-\alpha }\left\|
u_{0}\right\| _{p}\right) ^{q^{*}}+3^{q^{*}-1}\left( C_{11}\delta ^{-\alpha
}\right) ^{q^{*}}\int\limits_{0}^{t}h_{1}^{q^{*}}(s)W(s)^{lq^{*}}ds \\
+3^{q^{*}-1}\left( C_{12}\delta ^{-\alpha }\right)
^{q^{*}}\int\limits_{0}^{t}h_{2}^{q^{*}}(s)U(s)^{lq^{*}}ds,%
\end{array}
\label{13}
\end{equation}
and
\begin{equation}
\begin{array}{c}
W(t)^{q^{*}}\leq 5^{q^{*}-1}\left( C_{3}\delta ^{-\alpha }\left\|
w_{0}\right\| _{p}\right)
^{q^{*}}+5^{q^{*}-1}C_{10}^{q^{*}}\int%
\limits_{0}^{t}h_{4}^{q^{*}}(s)W(s)^{q^{*}}ds \\
+5^{q^{*}-1}\left( C_{13}\delta ^{-\alpha }\right)
^{q^{*}}\int\limits_{0}^{t}h_{1}^{q^{*}}(s)W(s)^{lq^{*}}ds+5^{q^{*}-1}\left(
C_{14}\delta ^{-\alpha }\right)
^{q^{*}}\int\limits_{0}^{t}h_{3}^{q^{*}}(s)W(s)^{lq^{*}}ds \\
+5^{q^{*}-1}\left( C_{15}\delta ^{-\alpha }\right)
^{q^{*}}\int\limits_{0}^{t}h_{2}^{q^{*}}(s)U(s)^{lq^{*}}ds.%
\end{array}
\label{14}
\end{equation}
Putting $F(t)=U(t)^{q^{*}}+W(t)^{q^{*}},$ we infer from \eqref{13} and \eqref{14} that
\begin{equation*}
F(t)\leq C_{0}\left( \left\| u_{0}\right\| _{p},\left\| w_{0}\right\|
_{p}\right) +\int\limits_{0}^{t}h(s)(F(s)+F(s)^{l})ds,\;t\geq \delta >0
\end{equation*}
where $C_{0}\left( \left\| u_{0}\right\| _{p},\left\| w_{0}\right\|
_{p}\right) =3^{q^{*}-1}\left( C_{1}\delta ^{-\alpha }\left\| u_{0}\right\|
_{p}\right) ^{q^{*}}+5^{q^{*}-1}\left( C_{3}\delta ^{-\alpha }\left\|
w_{0}\right\| _{p}\right) ^{q^{*}}$ and
$$
h(s)=\max \left\{
(C_{17}+C_{19})h_{1}^{q^{*}}(s),C_{20}h_{3}^{q^{*}}(s),(C_{18}+C_{21})h_{2}^{q^{*}}(s),C_{22}h_{4}^{q^{*}}(s)\right\}
$$
where $C_{i}$, $i=17,18,...,22$ are the coefficients of the integral terms
in \eqref{13} and \eqref{14} in the order.

Let $G(z)=\int\limits_{z_{0}}^{z}\frac{dy}{y+y^{l}}.$ Then, by Lemma \ref{L6} we
may conclude that
\begin{equation*}
\begin{array}{c}
F(t)\leq G^{-1}\left[ G\left( C_{0}\left( \left\| u_{0}\right\| _{p},\left\|
w_{0}\right\| _{p}\right) \right) +\int\limits_{0}^{t}h(s)ds\right] \\
\leq C_{0}\left( 1+C_{0}^{l-1}\right) ^{\frac{1}{1-l}}e^{\int%
\limits_{0}^{t}h(s)ds}.\left[ 1-\frac{C_{0}}{1+C_{0}^{l-1}}%
e^{(l-1)\int\limits_{0}^{t}h(s)ds}\right] ^{\frac{1}{1-l}}.%
\end{array}
\end{equation*}
From our assumptions on $\left\| u_{0}\right\| _{p},$ $\left\| w_{0}\right\|
_{p}$ and $h(t)$ we deduce that
\begin{equation}
\left\| A_{p}^{\alpha }u\right\| _{p}\leq Ce^{-(b-\varepsilon )t}\left(
\left\| u_{0}\right\| _{p}+\left\| w_{0}\right\| _{p}\right)  \label{15}
\end{equation}
and
\begin{equation}
\left\| G_{p}^{\alpha }w\right\| _{p}\leq Ce^{-(b-\varepsilon )t}\left(
\left\| u_{0}\right\| _{p}+\left\| w_{0}\right\| _{p}\right) .  \label{16}
\end{equation}
The decay rates in $C^{\mu }(\bar{\Omega}),$ $\mu \in [0,2)$ follow from
Lemma \ref{L4}.

\underline{B. The decay rate of $\left\| v-v_{\infty }\right\| _{C^{\mu }(%
\bar{\Omega})}:$}

As in Hoshino \cite{3}, let us write
\begin{equation*}
v-v_{\infty }=\left( Q_{0}v(t)-v_{\infty }\right) +Q_{+}v(t),
\end{equation*}
and estimate the terms in the right hand side separately.

\underline{a. The estimation of $Q_{0}v(t)-v_{\infty }:$}

Integrating the second equation in \eqref{1} over $(0,t)\times \Omega $, we have
\begin{equation*}
\int\limits_{\Omega }v(x,t)dx+\int\limits_{0}^{t}\int\limits_{\Omega
}\left\{ (a_{1}(s)+a_{3}(s))w^{m}-a_{2}(s)u^{n}v^{k}\right\}
dxds=\int\limits_{\Omega }v_{0}(x)dx.
\end{equation*}
It appears then that
\begin{equation*}
\left| Q_{0}v(t)-v_{\infty }\right| =\left| \Omega \right| ^{-1}\left|
\int\limits_{t}^{\infty }\int\limits_{\Omega }\left\{
(a_{1}(s)+a_{3}(s))w^{m}-a_{2}(s)u^{n}v^{k}\right\} dxds\right| .
\end{equation*}

In the rest of the proof, $C$ will denote a generic positive constant which
may be different at different occurrences.

Using \eqref{9}, \eqref{15}, \eqref{16} and the H\"{o}lder inequality we see that
\begin{equation}
\begin{array}{l}
\left| Q_{0}v(t)-v_{\infty }\right| \\
\leq C\int\limits_{t}^{\infty }\left\{ (s^{\sigma _{1}}h_{1}+s^{\sigma
_{3}}h_{3})(s)e^{-l(b-\varepsilon )s}+s^{\sigma
_{2}}h_{2}(s)M^{k}e^{-l(b-\varepsilon )s}\right\} ds \\
\leq C\int\limits_{t}^{\infty }\left\{ \sum\limits_{i=1}^{3}s^{\sigma
_{i}}h_{i}(s)\right\} e^{-l(b-\varepsilon )s}ds \\
\leq C\int\limits_{t}^{\infty }\left\{ \sum\limits_{i=1}^{3}s^{\sigma
_{i}}h_{i}(s)\right\} e^{-(b-\varepsilon )s}.e^{-(l-1)(b-\varepsilon )s}ds
\\
\leq Ce^{-(b-\varepsilon )t}\int\limits_{t}^{\infty }\left\{
\sum\limits_{i=1}^{3}s^{\sigma _{i}}h_{i}(s)\right\} e^{-(l-1)(b-\varepsilon
)s}ds \\
\leq Ce^{-(b-\varepsilon )t}.%
\end{array}
\label{17}
\end{equation}
We also used Lemma 7 in the last inequality.

\underline{b. The estimation of $Q_{+}v(t):$}

In order to estimate $Q_{+}v(t)$ let us apply $B_{p+}^{\alpha }Q_{+}$ to the
integral equation associated with the second equation in \eqref{1}. We find for
all $t\geq \delta >0$%
\begin{equation*}
\begin{array}{c}
B_{p+}^{\alpha }Q_{+}v(t)=B_{p+}^{\alpha }e^{-(t-\delta
)B_{p+}}Q_{+}v(\delta ) \\
+\int\limits_{\delta }^{t}B_{p+}^{\alpha }e^{-(t-s)B_{p+}}Q_{+}\left[
(a_{1}(s)+a_{3}(s))w^{m}-a_{2}(s)u^{n}v^{k}\right] ds.%
\end{array}
\end{equation*}
Taking the $L^{p}$-norm and using the second inequality in Lemma 5, we
obtain
\begin{equation*}
\begin{array}{c}
\left\| B_{p+}^{\alpha }Q_{+}v(t)\right\| _{p}\leq C_{2}(t-s)^{-\alpha
}e^{-d_{2}(t-s)\lambda }\left\| Q_{+}v(\delta )\right\| _{p}+C_{2}\left\|
Q_{+}\right\| _{L^{p}(\Omega )\rightarrow L^{p}(\Omega )} \\
\int\limits_{\delta }^{t}(t-s)^{-\alpha }e^{-d_{2}\lambda (t-s)}\left\{
(s^{\sigma _{1}}h_{1}+s^{\sigma _{3}}h_{3})\left\| w^{m}\right\| +s^{\sigma
_{2}}h_{2}M^{k}\left\| u^{n}\right\| \right\} ds.%
\end{array}
\end{equation*}
Next, using \eqref{9}, \eqref{15} and \eqref{16} we see that for all $t\geq \delta +T$%
\begin{equation}
\begin{array}{ll}
\left\| B_{p+}^{\alpha }Q_{+}v(t)\right\| _{p} & \leq Ce^{-d_{2}\lambda
t}\left\| Q_{+}v(\delta )\right\| _{p} \\
& +Ce^{-d_{2}\lambda t}\int\limits_{\delta }^{t}(t-s)^{-\alpha
}e^{d_{2}\lambda s}e^{-l(b-\varepsilon )s}\left(
\sum\limits_{i=1}^{3}s^{\sigma _{i}}h_{i}(s)\right) ds \\
& \leq Ce^{-d_{2}\lambda t}\left\{ \left\| Q_{+}v(\delta )\right\|
_{p}+\int\limits_{0}^{t-\delta }(t-\delta -s)^{-\alpha }\right. \\
& \left. \times e^{-\left[ l(b-\varepsilon )-d_{2}\lambda \right] (s+\delta
)}\left( \sum\limits_{i=1}^{3}(s+\delta )^{\sigma _{i}}h_{i}(s+\delta
)\right) ds\right\}.
\end{array}
\label{18}
\end{equation}
\begin{itemize}
\item[(i)] If $d_{2}\lambda <lb$, then choose $\varepsilon $ such that $%
0<l\varepsilon <lb-d_{2}\lambda .$ Hence, we may apply Lemma 8, together with
the H\"{o}lder inequality to get
\begin{equation}
\left\| B_{p+}^{\alpha }Q_{+}v(t)\right\| _{p}\leq Ce^{-d_{2}\lambda
t}\;for\;all\;t\geq \delta +T.  \label{19}
\end{equation}
\item[(ii)] If $d_{2}\lambda \geq lb,$ then $l(b-\varepsilon )-d_{2}\lambda <0.$
Multiplying by $e^{\rho (s+\delta )}.e^{-\rho (s+\delta )},$ with $\rho
>d_{2}\lambda -l(b-\varepsilon ),$ the integrand in \eqref{18}, using the H\"{o}lder inequality and Lemma 8, we see that
\begin{eqnarray}
\left\| B_{p+}^{\alpha }Q_{+}v(t)\right\| _{p} &\leq &Ce^{-d_{2}\lambda
t}\left\{ \left\| Q_{+}v(\delta )\right\| _{p}+e^{-\tilde{\rho}t}\right\}
\label{20} \\
&\leq &Ce^{-(d_{2}\lambda -\tilde{\rho})t}.  \notag
\end{eqnarray}
The conclusion follows from \eqref{17}, \eqref{19} and \eqref{20}.
\end{itemize}
\end{proof}
\section{Concluding Remarks}
Some remarks are in order:
\begin{itemize}
\item The condition $a_{1}(t)\leq Ca_{3}(t)$ needed in Theorem \ref{T1} has not been
used in Theorem \ref{T2}. So, the decay rates hold provided one may prove a
convergence result without this condition.
\item The assumption $\max \left\{ 1,\frac{N(m-1)}{2p\alpha },\frac{N(n-1)}{2p\alpha }\right\} <\min \{m,n\}$ may be relaxed somewhat using different $l_{1}$ and $l_{2}$ and applying Lemma 9 with $p=2.$
\item We also have exponential decay in (b) (ii) without the growth condition on $h_{i}$, $i=1,2,3$ in case $\sigma _{i}=0,$ $i=1,2,3.$
\item It is possible to obtain sharper estimates using comparison results by
replacing the bounds $M^{k}$ with $v_{\infty }^{k}-Ce^{-(b-\varepsilon )t}$
in \eqref{5} and $v_{\infty }^{k}+Ce^{-(b-\varepsilon )t}$ in \eqref{6}, for large
values of $t$.
\end{itemize}

\section*{Acknowledgment}
We express our deep thanks to {\em Bao Quoc Tang} for giving us the proof of Theorem \ref{GE} and for interesting discussions. The second author would like to thank King Fahd
University of Petroleum and Minerals (Saudi Arabia) for its support.

\end{document}